\documentclass[review]{elsarticle}\usepackage[utf8]{inputenc}
\usepackage[margin=1in]{geometry}
\usepackage[hidelinks]{hyperref}
\usepackage{amsfonts}
\usepackage{amsmath}
\usepackage{mathrsfs} 
\usepackage{amssymb,amsmath}

\usepackage{mathptmx}       % selects Times Roman as basic font
\usepackage{helvet}         % selects Helvetica as sans-serif font
\usepackage{courier}        % selects Courier as typewriter font
\usepackage{type1cm}        % activate if the above 3 fonts are
\usepackage{makeidx}         % allows index generation
\usepackage{graphicx}        % standard LaTeX graphics tool
\usepackage{multicol}        % used for the two-column index
\usepackage[bottom]{footmisc}% places footnotes at page bottom
\usepackage{subfigure}
\usepackage{booktabs}

\usepackage{epsfig} % for postscript graphics files
\usepackage{epstopdf}
\usepackage{color}
\usepackage{colordvi}
\usepackage[dvipsnames]{xcolor}
\usepackage{lineno}
\newcommand{\bfm}[1]{ \mathbf{#1}     }     
\newcommand{\Nabla}{ \boldsymbol{\nabla} }   % Vector Nabla

\newcommand{\dx}{ \; {\rm d x}    }                % dx (vector)
\newcommand{\dy}{ \, {\rm d y}    }                          % dx (vector)
\newcommand{\dX}{ \; {\rm d } \bfm{x}   }  
\newcommand{\dt}{ \, {\rm d t}    }
\newcommand{\dss}{ \, {\rm d s}    }
\newcommand{\VV}{\mathcal{V}(\Omega)}
\newcommand{\UU}{\mathcal{U}(\Omega)}
\newcommand{\UUT}{\mathcal{U}(\Omega_T)}
\newcommand{\VVT}{\mathcal{V}(\Omega_T)}
\newcommand{\VVh}{\mathcal{V}_h(\Omega)}
\newcommand{\UUh}{\mathcal{U}_h(\Omega)}
\newcommand{\ds}{\displaystyle}

\newcommand{\eirik}[1]{\textcolor{black}{#1}}

\newcommand{\matt}[1]{\textcolor{black}{#1}}
\newcommand{\one}[1]{\textcolor{black}{#1}}
\newcommand{\two}[1]{\textcolor{black}{#1}}

\begin{document}

\begin{frontmatter}
 \title{Extending FEniCS to Work in Higher Dimensions Using Tensor Product Finite Elements}

\author[1]{Mark Loveland\corref{cor1}}
\ead{Markloveland@utexas.edu}

\author[1]{Eirik Valseth}
%\ead{Eirik@utexas.edu}

\author[2]{Matt Lukac}
%\ead{mlukac@uoregon.edu}

 \author[1]{Clint Dawson }
%\ead{Clint.Dawson@austin.utexas.edu}

%\tnotetext[t2]{This work has been supported by the..........   }

\cortext[cor1]{Corresponding author}

 %\fntext[fn1]{This is the first author footnote.}
%\ \fntext[fn2]{Another author footnote, this is a very long 
%   footnote and it should be a really long footnote. But this 
%   footnote is not yet sufficiently long enough to make two 
%   lines of footnote text.}

 \address[1]{Oden Institute for Computational Engineering and Sciences, The University of Texas at Austin. 201 E. 24th St. Stop C0200, Austin, Texas 78712  }
  \address[2]{Institute of Ecology and Evolution, Department of Biology, University of Oregon. 77 Klamath Hall, 1210 University of Oregon, Eugene, Oregon 97403}

\begin{keyword}
 FEniCS, tensor product, Cartesian product, finite element method \MSC 65N30 ,65M60, 35D30
\end{keyword}

\biboptions{sort&compress}

%\newpageafter{title} 

%
%\maketitle
%
%
\begin{abstract}
We present a method to extend the finite element library FEniCS to solve problems with domains in \eirik{dimensions above three by constructing} tensor product finite elements. This methodology only requires that the high dimensional domain is structured as a Cartesian product of two lower dimensional subdomains. In this study we consider \one{Dirichlet problems for scalar} linear partial differential equations, though the methodology can be extended to non-linear problems. The utilization of tensor product finite elements allows us to construct a global system of \eirik{linear algebraic} equations that only \eirik{relies on} the finite element infrastructure of the lower dimensional subdomains contained in FEniCS. 
\eirik{We demonstrate the effectiveness of our methodology} in \one{four} \eirik{distinctive} test cases. The first test case is a Poisson equation  \eirik{posed in} a four dimensional domain which is a Cartesian product of two unit squares \eirik{solved using the classical Galerkin finite element method}. The second test case is the wave equation in space-time, where the \eirik{computational} domain is a Cartesian product of a two dimensional space grid and a one dimensional time interval. \eirik{In this second case we also employ the Galerkin method. } The third test case is an advection dominated advection-diffusion equation where the global domain is a Cartesian product of two one dimensional intervals  in which the streamline upwind Petrov-Galerkin method \eirik{is applied to ensure discrete stability}. \one{The final test case uses the Galerkin approach to solve a Poisson problem on a Cartesian product of two intervals with a spatially varying, non-separable diffusivity term.} In all cases, \two{a p=1 basis is used and} optimal $L^2$ convergence rates of order $h^{p+1}$ \two{of the errors } are achieved with respect to $h$ refinement. 
\end{abstract}

\end{frontmatter}

%\maketitle

\section{Introduction}

In the last decade, open source libraries with high level APIs that automate the process of solving partial differential equations (PDEs) using the finite element method (FEM) have become valuable tools in computational research. Some prominent libraries of this kind include Firedrake~\cite{rathgeber2016firedrake}, deal.II~\cite{bangerth2007deal}, MFEM~\cite{anderson2021mfem}, and  FEniCS~\cite{alnaes2015fenics}. FEniCS is one of the most widely used of these libraries. FEniCS, along with most other FEM libraries, only include capability of handling up to three dimensional problems (GetFEM~\cite{renard2020getfem} is a notable exception). FEniCS specifically supports unstructured meshes with corresponding basis functions in one, two, and three dimensions. Extension of FEniCS to solve problems in higher dimensions is of interest for problems such as spectral wind wave models, spatial population genetics, quantum mechanics, and even 3d space-time problems.

FEniCS does not contain support for fully unstructured meshes in dimensions higher than three. However, at the expense of losing a fully unstructured mesh, a higher dimensional space can be discretized with the available tools in FEniCS if it is the Cartesian product of two lower dimensional (3 or lower) spaces that can themselves be unstructured. It is well known that in a Cartesian product space, a finite element basis can be constructed called the product basis which spans the function space of the full domain~\cite{brenner2008mathematical,quarteroni2010numerical,ern2013theory}. This research seeks to use the infrastructure of an FEM library, such as FEniCS, in the lower dimensional spaces to construct a basis for this high dimensional space. 

 Using FEM to discretize a structured domain via a product basis has been done before without the use of FEM libraries such as FEniCS. In 1978, Banks~\cite{Bank1978} used tensor product finite elements to solve a 2D Poisson problem on a structured grid in order to find a faster solver. In the study, the domain was a Cartesian product between two unit intervals resulting in a uniform square mesh. Banks used the tensor product of two one dimensional quadratic and cubic basis functions to construct the polynomial basis of the full 2D domain. In 1980, Baker performed numerical tests as well as a stability analysis on a tensor product finite method with application to convection-dominated fluid flow problems~\cite{BAKER1981215}. The stability analysis showed the basic algorithm is spatially fourth- order accurate in its most elementary embodiment and the numerical experiments on the convection-dominated model test problems confirmed the basic viability of the developed algorithm, and its tensor product formulation. Recently, Du \eirik{\emph{et al.}~\cite{DU2013181}} used tensor product finite elements to construct a fast solver for an electromagnetics scattering problem of a cavity. In  Firedrake~\cite{rathgeber2016firedrake}, there exists a capability to create tensor product finite elements only up to three dimensions~\cite{mcrae2016automated}. A package built on deal.II called Hyperdeal~\cite{munch2020hyperdeal} has the capability of creating tensor product finite elements in up to six dimensions. deal.II is different to a library such as FEniCS since deal.II uses quadrilateral elements in two dimensions whereas FEniCS uses triangles. \eirik{Here, we present a general software framework built on the components of the FEniCS library as well as other open source Python libraries to construct high dimensional meshes and corresponding FE discretizations.  }
 %In this study, a FEM library that includes triangular elements in two dimensions such as FEniCS is considered.

Following this introduction, the general set of problems for which this paper focuses on will be defined and the notation for the product basis will be introduced \eirik{in Section~\ref{sec:model_prob}}. Then, a set of four different model problems will be described in detail and the derivation of the system of equations using tensor product finite elements will be shown for each case. The first model problem will be discussed in Section~\ref{sec:model_prob} which is an N dimensional Poisson problem where the domain is a Cartesian product between two subdomains of dimension three or less. In Section~\ref{sec:prob2} the second model problem is discussed which is a wave equation where the domain is decomposed as a Cartesian product between space and time. The third model problem is discussed in Section~\ref{sec:prob3} which is an advection dominated advection-diffusion equation where the advection aligns with one of the subdomains, a streamlined upwind Petrov-Galerkin method (SUPG) is formulated for this case. In Section~\ref{sec:prob4} the fourth model problem is discussed which is a Poisson problem in 2d which has a spatially varying diffusivity term that is non-separable. After each system of equations is derived, numerical tests were run using FEniCS for each model problem where specific boundary conditions were given. For each case, error and convergence rates are tabulated in Section~\ref{sec:numerical_results}. Lastly, conclusions and recommendations for future work are given in Section~\ref{sec:Conclusions}.
% (list what we present in this paper in which section, see for example one of my papers)

\section{Methods} \label{sec:methods}

\begin{figure}[h!]
    \centering
    \includegraphics[width=\textwidth]{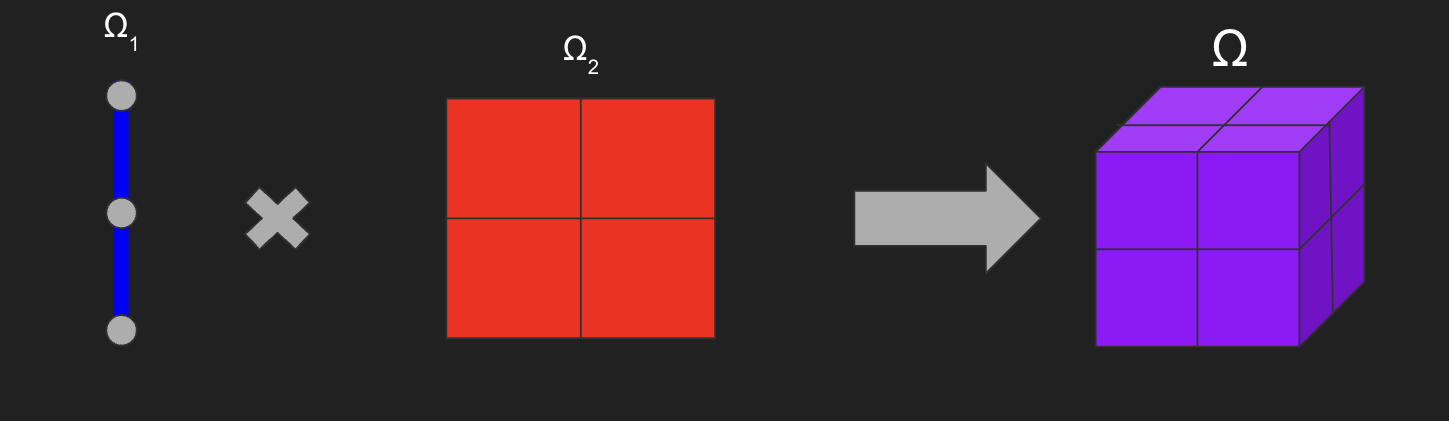}
    \caption{Example of a global domain $\Omega$ that is a Cartesian product of two lower dimensional subdomains $\Omega_1, \Omega_2$.}
    \label{fig:cartesian_illustration}
\end{figure}

To present the proposed methodology and algorithms, we 
consider the following class of problems, i.e., PDEs:
\begin{equation*}
\begin{split}
    \mathscr{L} u = f \quad \textrm{in} \quad \Omega = \Omega_1 \times \Omega_2.
\end{split}
\end{equation*}
 Where $\mathscr{L}$ is a linear differential operator, $f$ is a forcing function, and the domain $\Omega$ is defined as a Cartesian product between two lower dimensional Lipschitz domains.  For example, if the global domain $\Omega  \subset \mathbb{R}^2$, then $\Omega$ can be defined by the Cartesian product of 2 intervals $\Omega_1\subset \mathbb{R}$ and  $\Omega_2 \subset \mathbb{R}$. Hence, we  consider general domains of the form $\Omega = \Omega_1 \times \Omega_2 $, see Figure~\ref{fig:cartesian_illustration} for an illustration.
FEniCS can discretize up to three dimensional objects, hence,  in practice this framework can be used to define domains that are Cartesian products of up to six dimensions. Furthermore, this process can be done iteratively to yield even higher dimensional domains as $\Omega_1$ itself could be a product of two other spaces and so on. We note that in this presentation, we consider only symmetric functional settings which admit well posed  Galerkin FE discrtetization. \one{For the methods described in model Sections~\ref{sec:model_prob},~\ref{sec:prob2}, and~\ref{sec:prob3}, $\mathscr{L}$ is restricted to linear differential operators that either do not contain functions which depend on domain coordinates or contain functions that are separable.   However, it is important to note that a generalization of these methods does apply to operators that contain non-separable functions. The procedure is a bit more complicated and some efficiency is lost since the use of Kronecker products to directly assemble the system as will be seen in Sections~\ref{sec:model_prob},~\ref{sec:prob2}, and~\ref{sec:prob3} is no longer possible but an example implementation is constructed in Section~\ref{sec:prob4}. }

It can be shown that if a function $u$ is defined on a Cartesian product domain $\Omega = \Omega_1 \times \Omega_2$ then we can construct a basis of a polynomial space to be used for the FEM \eirik{ by exploiting this underlying geometric structure, see, e.g.,}~\cite{quarteroni2010numerical,brenner2008mathematical,ern2013theory}. 
 \eirik{With} a basis for a polynomial space for the first lower dimensional subdomain $\Omega_1$:
\begin{equation}
    \{\phi_i\}_{i=1}^N,
\end{equation}
and a basis for the second subdomain $\Omega_2$:
\begin{equation}
    \{\psi_j\}_{j=1}^M.
\end{equation}
The basis for the entire domain can subsequently be constructed and any arbitrary function $u$ whose domain is in $\Omega = \Omega_1 \times \Omega_2$ can be approximated via the product basis:
\begin{equation}
    u \approx \sum_{i=1}^N\sum_{j=1}^M u_{i,j} \phi_i\, \psi_j. 
\end{equation}
\eirik{Note that in the following,} we use boldface letters and symbols to denote vector quantities, e.g., $\dX = \dx \dy$. 

\subsection{Model Problem 1: N dimensional Poisson Equation} \label{sec:model_prob}
As a first model problem to illustrate the methodology, we consider the Poisson equation:
\begin{equation}\label{eqn:model_prob1}
\begin{split}
        - \Delta \, u &= f \quad \textrm{in} \quad \Omega, \\
        u &= u_D \quad \textrm{on} \quad \partial \Omega, 
\end{split}
\end{equation}
where the source $f$ is in $L^2(\Omega)$ and the source data $u_D$ is assumed to be sufficiently regular. ~\one{Now we will briefly derive the finite element formulation using the Galerkin approach and the product basis. Note, that this derivation is not new and similar derivations can be found in the literature such as the work from Banks~\cite{Bank1978} for example. The derivation is included so it is clear how to implement in an algorithm as well as how the ideas of using tensor product elements will apply to more complex cases.} 
To define the weak formulation for~\eqref{eqn:model_prob1}, multiply both sides with a test function $v$ in $L^2(\Omega)$ and integrate over $\Omega$:
\begin{equation}\label{eqn:weak_form1}
    \int_{\Omega}  - \Delta \, u\, v \dX = \int_{\Omega} f \, v \dX \quad \forall v \in L^2(\Omega), 
\end{equation}
Integrating by parts on the left side (assuming Dirichlet boundary conditions on entire boundary) gives the following weak formulation\matt{: find $u \in \UU$ such that}
\begin{equation} \label{eq:weak_form}
    \int_{\Omega}  \Nabla u \cdot \Nabla v \dX = \int_{\Omega} f\,v \dX \quad \forall \quad v \in \VV, 
\end{equation}
where the function space $\VV$ is the Hilbert space with zero trace on the boundary  $H^1_0(\Omega)$ and $\UU$ is $\VV$ with a finite energy lift on the boundary so that the Dirichlet condition $u=u_D$ is satisfied.
The weak formulation in~\eqref{eq:weak_form}, and its corresponding discretization is known to be well posed, see, e.g.,~\cite{becker1981finite}. 

With a well posed weak formulation at hand, we can discretize this weak form using a finite element basis. In this case, the functional setting dictates the use of a $C^0$ continuous \eirik{polynomial} basis for the classical FEM.
Hence, we approximate the trial functions with the product basis:
\begin{equation}
    u \approx \sum_{i=1}^N\sum_{j=1}^M u_{i,j} \, \phi_i \, \psi_j,
\end{equation}
and the test functions by its product basis:
\begin{equation}
    v \approx \sum_{k=1}^N\sum_{l=1}^M v_{k,l} \, \gamma_k \, \beta_l.
\end{equation}
The weak form from~\eqref{eqn:weak_form1} can then be discretized by substitution of the product bases:
\begin{equation}
\ds    \int_{\Omega}  (\Nabla \sum_{i=1}^N\sum_{j=1}^M u_{i,j} \, \phi_i \, \psi_j) \cdot ( \Nabla \sum_{k=1}^N\sum_{l=1}^M  v_{k,l}\,\gamma_k\, \beta_l) \dX = \int_{\Omega} f\, (\sum_{k=1}^N\sum_{l=1}^M  v_{k,l}\,\gamma_k\, \beta_l ) \dX.
\end{equation}
Due to arbitrariness of the test function and the bilinearity of the weak form, this implies the following form:
\begin{equation}
    \sum_{i=1}^N\sum_{j=1}^M \int_{\Omega}  \Nabla   \, (\phi_i \, \psi_j) \cdot  \Nabla \,(\gamma_k\, \beta_l) \dX \, u_{i,j} = \int_{\Omega} f\, (\,\gamma_k\, \beta_l ) \dX.
\end{equation}
This left hand side \eirik{results} in a \eirik{product between a} 4 dimensional \eirik{and} a 2 dimensional tensor $u_{i,j}$. \eirik{Consequently,} the 4 dimensional tensor $A_{ijkl}$ will be of the form:
\begin{equation}
        \int_{\Omega}   \Nabla  (\phi_i \,\psi_j) \cdot \Nabla  (\gamma_k \,\beta_l \,) \dX.
\end{equation}
By the product rule we have:
\begin{equation} \label{eq:chain_ruled}
        \int_{\Omega}   \Nabla  (\phi_i \, \psi_j) \cdot \Nabla  (\gamma_k \,\beta_l) \dX  = \int_{\Omega}\, (\phi_i \, \Nabla (\psi_j) + \psi_j\, \Nabla (\phi_i)) \cdot (\gamma_k\, \Nabla (\beta_l) + \beta_l\, \Nabla (\gamma_k) )\dX.
\end{equation}
By construction, \eirik{the} $\phi$\eirik{'s}  and$\gamma$\eirik{'s} only vary in the first subdomain $\Omega_1$ whereas \eirik{the} $\psi$\eirik{'s} and $\beta$\eirik{'s} only vary in the second subdomain $\Omega_2$. \eirik{Thus, the integral form~\eqref{eq:chain_ruled} can be simplified. To this end, we use the following notation convention:} if the gradient operator on the entire domain \eirik{$\Omega$} is $\Nabla$, define the gradient on the subdomains \eirik{$\Omega_1$ and $\Omega_2$} \eirik{as} $\Nabla_1$, $\Nabla_2$\eirik{, respectively}. \eirik{Hence, by} construction $\Nabla = (\Nabla_1, \Nabla_2)$.
Rewriting \eirik{ the gradients in~\eqref{eq:chain_ruled} then} gives:
\begin{equation} \label{eq:expansion}
\begin{split}
   & \int_{\Omega} (\phi_i \, (\Nabla_1, \Nabla_2)\,(\psi_j) + \psi_j \,(\Nabla_1, \Nabla_2)\, (\phi_i)) \cdot (\gamma_k\, (\Nabla_1, \Nabla_2)\, (\beta_l) + \beta_l (\Nabla_1, \Nabla_2)(\gamma_k))\dX,  \\ \qquad &=
    \int_{\Omega} (  (\mathbf{0}, \phi_i\,\Nabla_2 (\psi_j)) +  (\psi_j\,\Nabla_1(\phi_i), \mathbf{0}) ) \cdot ( (\mathbf{0}, \gamma_k\, \Nabla_2(\beta_l))  +  (\beta_l\,\Nabla_1(\gamma_k), \mathbf{0}))\dX,  \\ &=
    \int_{\Omega}  (\psi_j\,\Nabla_1(\phi_i), \phi_i\,\Nabla_2 (\psi_j))  \cdot  (\beta_l\,\Nabla_1(\gamma_k), \gamma_k\, \Nabla_2(\beta_l))  \dX, \\ \qquad &= 
    \int_{\Omega}  \Nabla_1(\phi_i) \cdot \Nabla_1(\gamma_k)\,\psi_j\,\beta_l +  \phi_i \,\gamma_k \,\Nabla_2 (\psi_j) \cdot \Nabla_2\,(\beta_l)  \dX. 
\end{split}
\end{equation}
%
%
%Multiplying out the dot product we get:
%
%
%\begin{equation}
%       \int_{\Omega}  \Nabla_1\phi_i \cdot \Nabla_1\gamma_k\psi_j\beta_l +  \phi_i \gamma_k \Nabla_2 \psi_j \cdot \Nabla_2\beta_l  \dX  
%\end{equation}
%
%
Since \eirik{the} domain $\Omega$ is a Cartesian product of \eirik{the} subdomains $\Omega_1, \Omega_2$ we can rewrite the \eirik{last} integral \eirik{in~\eqref{eq:expansion}, %\eirik{(note that $\dx$ and $\dy$ denote integration over $\Omega_1$ and $\Omega_2$, respectively)
as}:
\begin{equation}\label{eqn:product_rule}
           \int_{\Omega_1} \int_{\Omega_2}  \Nabla_1(\phi_i) \cdot \Nabla_1(\gamma_k\,)\psi_j\,\beta_l +  \phi_i\, \gamma_k\, \Nabla_2 (\psi_j) \cdot \Nabla_2(\beta_l)  \dy  \dx.
\end{equation}
\eirik{An application of} Fubini's theorem \eirik{gives}:
\begin{equation}
               \int_{\Omega_1} \Nabla_1(\phi_i) \cdot \Nabla_1(\gamma_k) \dx \int_{\Omega_2}  \psi_j\,\beta_l \dy + \int_{\Omega_1} \phi_i\, \gamma_k \dx \int_{\Omega_2}  \Nabla_2 (\psi_j) \cdot \Nabla_2(\beta_l)   \dy,
\end{equation}
These can be written as Kronecker products of  matrices,  \eirik{e.g., stiffness matrices computed with} FEniCS. \eirik{To continue the discussion, we introduce}  the following notation:
\begin{equation}
    \begin{split}
        K_{11} = \int_{\Omega_1} \Nabla_1\phi_i \cdot \Nabla_1\gamma_k \dx, \quad K_{22} = \int_{\Omega_2}  \psi_j\;\beta_l \dy,\\
        K_{12} = \int_{\Omega_1} \phi_i\; \gamma_k \dx, \quad
        K_{21} = \int_{\Omega_2}  \Nabla_2 \psi_j \cdot \Nabla_2\beta_l   \dy,
    \end{split}
\end{equation}
\eirik{where} the first index indicates subdomain (\eirik{i.e.,} 1 or 2 \eirik{in this case}) and the second indicates the term from the weak form (1 is the first order operator and 2  the \eirik{$0^{th}$} order). \eirik{Hence, each local matrix $K_{ij}$ can be constructed independently, and the global stiffness matrix $A_{ijkl}$ can be defined by the sum of 2 Kronecker products.}
%Since each of these matrices can be constructed then the global stiffness matrix $A_{ijkl}$ can be written as a sum of 2 kronecker products. 
The Kronecker product allows us to represent \eirik{the matrix of a} 4D tensor  \eirik{by smaller 2D matrices}:
\begin{equation}
    A = K_{11} \otimes K_{22} + K_{12} \otimes K_{21}.
\end{equation}
\one{For an arbitrary forcing function $f$, it is often more convenient to approximate $f$ as a member of the solution space.}  \eirik{A similar reasoning for the right hand side leads to the following}:
\one{
\begin{equation}
    \int_{\Omega}  f\;\gamma_k\; \beta_l \dX  \approx \sum_{i=1}^N\sum_{j=1}^M \int_{\Omega_1} \phi_i\; \gamma_k \dx \int_{\Omega_2} \psi_j\;\beta_l \dy, f_{i,j},
\end{equation}
}
\one{where $f_{i,j}$ is the pointwise value of $f$ at each $i,j$ coordinate in the global space. The forcing vector is consequently defined as:
\begin{equation}
    \mathbf{F} = f_{i,j}.
\end{equation}
}
\eirik{Finally,} the entire system \eirik{of equations} becomes:
\begin{equation} \label{eq:final_poisson_form}
    (K_{11} \otimes K_{22} + K_{12} \otimes K_{21})\mathbf{(U)} = K_{12} \otimes K_{22} \mathbf{F},
\end{equation}
where $\mathbf{U}$ are the values at all d.o.f in the Cartesian product space for $u$. 
\eirik{A subsequent application of boundary conditions to~\eqref{eq:final_poisson_form}  leads to the final system of linear algebraic equations.}

\subsection{Model Problem 2: Arbitrary dimensional Space-Time Wave Equation}\label{sec:prob2}
\eirik{The second model problem we consider is the linear wave equation. We consider this transient problem to highlight the application of this methodology to transient problems where space-time finite elements are employed. Thus, we have the following model problem:}
%Governing equation:
%
%
\begin{equation}\label{eqn:model_prob2}
\begin{split}
        \frac{\partial^2 u}{\partial t}-c^2\Delta u & = 0 \quad \textrm{in} \quad \Omega_T, \\
        \eirik{u} &\eirik{ =u_D \quad \textrm{on} \quad \partial \Omega,} \\
        \eirik{u} &\eirik{ =u_{initial} \quad \textrm{on} \quad \Omega,}
\end{split}
\end{equation}
\eirik{where $c$ denotes the wave speed}.
In this problem the space-time domain $\Omega_T$ is a Cartesian product domain of one, two, or three dimensional \eirik{spatial domain and a one dimensional temporal domain} ($\Omega_T = \Omega \times (0,T)$). \one{Now we will derive the finite element formulation for this problem, again using the Galerkin approach. Note that similar derivations can be found in the literature, see work of Loscher~\cite{Loscher2021} for instance.} The weak form of~\eqref{eqn:model_prob2} is obtained by multiplying by a test function and integrating over the entire space-time domain:
\begin{equation}
    \int_{\Omega_T} \frac{\partial^2 u}{\partial t} v  - c^2\Delta u v \dX \dt = 0 \quad \forall \quad v \in L^2(\Omega_T),
\end{equation}
\eirik{and subsequent integration } by parts in \eirik{space and time}  gives the following weak form. . Find $u \in \UUT$:
\begin{equation}\label{eqn:weak_form2}
   \int_{\Omega_T} -\frac{\partial u}{\partial t} \frac{\partial v}{\partial t} + c^2\Nabla u \cdot \Nabla v \dX\dt + \int_{\partial \Omega_T|_{t=T}}\frac{\partial u}{\partial t}v  \dX = 0 \quad \forall \quad v \in \VVT,
\end{equation}
where \eirik{we have applied Dirichlet conditions to the space-time boundary $\partial \Omega_T$, except at the final time boundary.} $\VVT$ is the space of all $H^1_0(\Omega_T)$ functions except on $\partial \Omega_T|_{t=T}$, where the trace is an unkown, and $\UUT$ is the space of all $H^1_0(\Omega_T)$ plus the trace $u=u_D$ on $\partial \Omega_T$ except the aforementioned part of the boundary $\Omega_T|_{t=T}$. 

\eirik{The corresponding discretization of}~\eqref{eqn:weak_form2} using a product basis gives a very similar system of equations to the Poisson problem \eirik{considered in Section~\ref{sec:model_prob}}, with a slight variation. The term $K_{21}$ must include \eirik{additional integrals} and becomes:
%
%Mark, will you check on the integration limits here ?
\begin{equation}
    K_{21} = \int_0^T  \frac{\partial \psi_j}{\partial t} \frac{\partial \beta_l}{\partial t}   dt - \int_{\partial \Omega_2} \Nabla \psi_j \cdot \mathbf{n} \beta_l \dss.
\end{equation}
\eirik{Consequently}, the global system of equations becomes:
\begin{equation} \label{eq:global_wave_sys}
    (c^2 K_{11} \otimes K_{22} - K_{12} \otimes K_{21})\mathbf{U} = \mathbf{0}.
\end{equation}
\eirik{Finally, application of boundary and initial conditions to~\eqref{eq:global_wave_sys} results in the final system of linear algebraic equations.} 
%Then the final system of equations is constructed by replacing rows in the system of equation corresponding to dirichlet boundaries with the appropriate row of the identity matrix and replacing the same row of the $\mathbf{0}$ vector with the value on the boundary.

\subsection{Model Problem 3: SUPG Stabilized Advection Dominated Advection Diffusion Equation}\label{sec:prob3}

\eirik{To highlight the versatility of our approach to consider non-standard FE techniques, we consider a PDE which is known to lead to stability issues in the Galerkin FE setting. Hence, we consider an advection-diffusion  PDE  in which advection is the dominant:}  
\begin{equation}\label{eqn:model_problem3}
\begin{split}
        -\kappa \, \Delta u + \mathbf{b} \cdot \Nabla u = f, \quad \textrm{on} \quad \Omega \\
         u = u_D \quad \textrm{on} \quad \partial \Omega,,
\end{split}
\end{equation}
where $\Omega$ is defined as a Cartesian Product of two lower dimensional Lipschitz domains: $\Omega = \Omega_1 \times \Omega_2$. \eirik{By following the standard procedure of deriving integral formulations,} we get the corresponding weak form. \one{Note that a similar derivation for a more complex advection problem can be found in the work of Baker~\cite{BAKER1981215}}:
%
%
%\begin{equation}
%    \int_{\Omega} -\kappa\, \Delta u\, v + \mathbf{b} \cdot \Nabla u\, v \dX = %\int_{\Omega}f\,v \dX \quad \forall v \in L^2(\Omega)
%\end{equation}
%
%
%Next, integrate the first term by parts, and assuming Dirichlet boundary conditions on the %entire domain boundary:
%
%
\begin{equation}\label{eqn:weak_form3}
        \int_{\Omega} \kappa \Nabla u \cdot \Nabla v + \mathbf{b} \cdot \Nabla (u) v dx = \int_{\Omega}fv dx \quad \forall v \in \VV
\end{equation}
In problems where the advection term dominates the diffusion, that is when Peclet number $Pe = \frac{L \| \mathbf{b} \| }{\kappa} >> 1$, where $L$ is the characteristic length, the standard Galerkin \eirik{method applied to~\eqref{eqn:weak_form3} may result in a}   discretization that is \eirik{unstable. } \eirik{This issue of stability can be overcome by careful design of the FE mesh or through stabilization techniques that ensure satisfaction of the discrete \emph{inf-sup} condition. 
Here, we consider the SUPG method introduced by Brooks and Hughes~\cite{brooks1982streamline} since it is widely used and has well developed criteria for discrete stability.}
\eirik{The SUPG method leads to stable FE discretizations by adjusting the discretized weak form~\eqref{eqn:weak_form3} with a penalized residual, i.e., find $u_h \in \UUh$: }
%To correct for his we add  the SUPG stability term so that (\ref{eqn:weak_form3}) becomes:
%
%
\begin{equation}\label{eqn:weak_form3b}
   \kappa(\Nabla u_h , \Nabla v_h)_{\Omega} + (\mathbf{b} \cdot \Nabla u_h, v_h)_{\Omega} + (\underbrace{-\kappa\,\Delta u_h + \mathbf{b} \cdot \Nabla u_h -f}_{\eirik{Residual}}, \tau(\mathbf{b} \cdot \Nabla v_h)  )_{\Omega} = (f,v_h)_{\Omega} \quad \forall v_h \in \VVh,
\end{equation}
\eirik{where the trial and test spaces consist of standard piecewise polynomials. Note that when the residual is zero, the stabilization term vanishes, i.e., it is consistent with the weak form~\eqref{eqn:weak_form3}.}

\eirik{As for the preceding model problems,} we wish to construct the 4 dimensional tensor $A_{ijkl}$ using the product bases of the test and trial spaces. \eirik{Substitution of} the product basis into (\ref{eqn:weak_form3b}) gives:
\one{
\begin{equation}
\begin{array}{c}
        \{\kappa(\Nabla (\phi_i\,\psi_j) , \Nabla (\gamma_k\,\beta_l))_{\Omega} + (\mathbf{b} \cdot \Nabla (\phi_i\, \psi_j), \gamma_k \beta_l)_{\Omega}
        + (-\kappa\Delta (\phi_i\,\psi_j) + \mathbf{b} \cdot \Nabla (\phi_i\,\psi_j), \tau[\mathbf{b} \cdot \Nabla(\gamma_k\,\beta_l)]  )_{\Omega}\} u_{i,j} = \\ (f,\,\gamma_k\, \beta_l + \tau[\mathbf{b} \cdot \Nabla(\gamma_k\,\beta_l)])_{\Omega}.
\end{array}
\end{equation}
}
\one{
\eirik{Application of} the product rule gives:
\begin{equation}
\begin{array}{c}
         \{\kappa(\psi_j \Nabla \phi_i + \phi_i \Nabla \psi_j , \beta_l \Nabla \gamma_k + \gamma_k \Nabla \beta_l)_{\Omega} + (\mathbf{b} \cdot (\psi_j \Nabla \phi_i +\phi_i \Nabla \psi_j), \gamma_k \beta_l)_{\Omega} +\\
         (-\kappa(\psi_j \Delta (\phi_i) + \phi_i \Delta (\psi_j)) + \mathbf{b} \cdot (\psi_j \Nabla (\phi_i) + \phi_i \Nabla(\psi_j)),
         \tau[\mathbf{b} \cdot (\beta_l \Nabla (\gamma_k) + \gamma_k \Nabla (\beta_l))]  )_{\Omega}\}u_{i,j} = \\ (f, \gamma_k  \beta_l + \tau[\mathbf{b} \cdot (\beta_l \Nabla (\gamma_k) + \gamma_k \Nabla (\beta_l))])_{\Omega}, 
\end{array}
\end{equation}
}
\one{\eirik{which we expand} using Fubini's theorem and approximate the exact f as a function in the discrete solution space:
\begin{equation} \label{eq:fubinisplit}
\begin{array}{l}
    \{\kappa \left[ (\Nabla \phi_i, \Nabla \gamma_k)_{\Omega_1}(\psi_j,\beta_l)_{\Omega_2} +  (\phi_i,\gamma_k)_{\Omega_1}( \Nabla \psi_j ,  \Nabla \beta_l)_{\Omega_2}\right] + (\mathbf{b} \cdot \Nabla \phi_i, \gamma_k)_{\Omega_1} (\psi_j,\beta_l)_{\Omega_2}
    + (\phi_i,\gamma_k)_{\Omega_1}(\mathbf{b} \cdot \Nabla \psi_j,\beta_l)_{\Omega_2} - \\
    \kappa \tau [(\Delta \phi_i, \mathbf{b} \cdot \Nabla \gamma_k)_{\Omega_1} (\psi_j \beta_l)_{\Omega_2}
    + (\Delta \phi_i , \gamma_k)_{\Omega_1}(\psi_j, \mathbf{b} \cdot  \Nabla \beta_l)_{\Omega_2} +    (\phi_i, \mathbf{b} \cdot \Nabla \gamma_k)_{\Omega_1}(\Delta \psi_j, \beta_l)_{\Omega_2} + (\phi_i,\gamma_k)_{\Omega_1}(\Delta \psi_j, \mathbf{b} \cdot \Nabla \beta_l)_{\Omega_2}] + \\
    \tau [ (\mathbf{b} \cdot \Nabla \phi_i, \mathbf{b} \cdot \Nabla \gamma_k)_{\Omega_1}(\psi_j,\beta_l)_{\Omega_2} + (\mathbf{b} \cdot \Nabla \phi_i, \gamma_k)_{\Omega_1}(\psi_j , \mathbf{b} \cdot \Nabla \beta_l)_{\Omega_2} + \\
    (\phi_i, \mathbf{b} \cdot \Nabla \gamma_k)_{\Omega_1}(\mathbf{b} \cdot \Nabla \psi_j, \beta_l )_{\Omega_2} + (\phi_i,\gamma_k)_{\Omega_1}(\mathbf{b} \cdot \Nabla \psi_j, \mathbf{b} \cdot \Nabla \beta_l)_{\Omega_2}]\}u_{i,j} \\
    = (\phi_i,\gamma_k)_{\Omega_1}(\psi_j,\beta_l)_{\Omega_2}f_{i,j} + \tau (\phi_i, \mathbf{b} \cdot \gamma_k)_{\Omega_1}(\psi_j, \beta_l)_{\Omega_2}f_{i,j} + \tau (\phi_i,\gamma_k)_{\Omega_1}(\psi_j,\mathbf{b} \cdot \Nabla \beta_l)_{\Omega_2}f_{i,j}.
\end{array}
\end{equation}
} %where $f_{i,j}$ represents the function value $f$ at each corresponding degree of freedom. }
The discrete weak form in~\eqref{eq:fubinisplit} can be represented as a global stiffness matrix compromised of the following smaller submatrices:
\one{
\begin{equation}
    \begin{split}
        K_{11} = (\Nabla \phi_i, \Nabla \gamma_k)_{\Omega_1} \quad  &
        K_{21} = ( \Nabla \psi_j ,  \Nabla \beta_l)_{\Omega_2}\\
        K_{12} = (\phi_i,\gamma_k)_{\Omega_1} \quad  &
        K_{22} = (\psi_j, \beta_l)_{\Omega_2}\\
        K_{13} = (\mathbf{b} \cdot \Nabla \phi_i,\gamma_k)_{\Omega_1}\quad  &
        K_{23} = (\mathbf{b} \cdot \Nabla \psi_j,\beta_l)_{\Omega_2}\\
        K_{14} = (\Delta \phi_i, \mathbf{b} \cdot \Nabla \gamma_k)_{\Omega_1}\quad  &
        K_{24} = (\Delta \psi_j, \mathbf{b} \cdot \Nabla \beta_l)_{\Omega_2}\\
        K_{15} = (\Delta \phi_i , \gamma_k)_{\Omega_1} \quad  &
        K_{25} = (\Delta \psi_j , \beta_l)_{\Omega_2}  \\
        K_{16} = (\phi_i, \mathbf{b} \cdot  \Nabla \gamma_k)_{\Omega_1}\quad  &
        K_{26} = (\psi_j, \mathbf{b} \cdot  \Nabla \beta_l)_{\Omega_2}\\
        K_{17} = (\mathbf{b} \cdot \Nabla \phi_i, \mathbf{b} \cdot \Nabla \gamma_k)_{\Omega_1}\quad  &
        K_{27} = (\mathbf{b} \cdot \Nabla \psi_j, \mathbf{b} \cdot \Nabla \beta_l)_{\Omega_2}\\
        F = f_{i,j} & \\
    \end{split}
\end{equation}
}
Then the global system of \one{linear algebraic equations is:} %equations to solve would be:
\one{
\begin{equation}
\begin{split}
    [ \kappa (K_{11} \otimes K_{22} + K_{12}\otimes K_{21}) + K_{13} \otimes K_{22} + K_{12} \otimes K_{23} - \kappa \tau(K_{14} \otimes K_{22} + K_{15} \otimes K_{26} + K_{16} \otimes K_{25} + K_{12} \otimes K_{24} ) + \\
    \tau ( K_{17} \otimes K_{22} + K_{13} \otimes K_{26} + K_{16} \otimes K_{23} + K_{12} \otimes K_{27} )]\mathbf{U} = (K_{12} \otimes K_{22} + \tau (K_{16} \otimes K_{22} + K_{12} \otimes K_{26}))\mathbf{F}.
\end{split}
\end{equation}
}

One advantage of this method is that it is possible to align one subdomain with the velocity vector. \eirik{In this special case} when $\mathbf{b}$ is only non-zero along one subdomain (e.g., $\Omega_2$) the above weak form~\eqref{eq:fubinisplit} reduces to:
\one{
\begin{equation}\label{eqn:weak_form3c}
\begin{split}
    \{\kappa \left[ (\Nabla \phi_i, \Nabla \gamma_k)_{\Omega_1}(\psi_j,\beta_l)_{\Omega_2} +  (\phi_i,\gamma_k)_{\Omega_1}( \Nabla \psi_j ,  \Nabla \beta_l)_{\Omega_2}\right]+ (\phi_i,\gamma_k)_{\Omega_1}(\mathbf{b} \cdot \Nabla \psi_j,\beta_l)_{\Omega_2} - \\
    \kappa \tau\left[ (\Delta \phi_i , \gamma_k)_{\Omega_1}(\psi_j, \mathbf{b} \cdot  \Nabla \beta_l)_{\Omega_2} + (\phi_i,\gamma_k)_{\Omega_1}(\Delta \psi_j, \mathbf{b} \cdot \Nabla \beta_l)_{\Omega_2}\right] \\
    + \tau \left[   (\phi_i,\gamma_k)_{\Omega_1}(\mathbf{b} \cdot \Nabla \psi_j, \mathbf{b} \cdot \Nabla \beta_l)_{\Omega_2}\right] \}u_{i,j} = \{(\phi_i,\gamma_k)_{\Omega_1}(\psi_j,\beta_l)_{\Omega_2} + \tau (\phi_i,\gamma_k)_{\Omega_1}(\psi_j,\mathbf{b} \cdot \Nabla \beta_l)_{\Omega_2}\}f_{i,j}.
\end{split}
\end{equation}
}
%Now we can construct a system of equations and solve for $u_{ij} $, we define the smaller stiffness matrices. First index represents the subdomain, second index represents the weak form term:
\eirik{Hence, the local matrices are defined:}
\begin{equation}
\begin{array}{cc}
        K_{11} = (\Nabla \phi_i, \Nabla \gamma_k)_{\Omega_1} \quad  &
        K_{21} = ( \Nabla \psi_j ,  \Nabla w)_{\Omega_2}\\
        K_{12} = (\phi_i,\gamma_k)_{\Omega_1}\quad  &
        K_{22} = (\psi_j, \beta_l)_{\Omega_2}\\
         K_{13} = (\Delta \phi_i , \gamma_k)_{\Omega_1} \quad  &
        K_{23} = (\mathbf{b} \cdot \Nabla \psi_j,\beta_l)_{\Omega_2}\\
        K_{24} = (\psi_j, \mathbf{b} \cdot  \Nabla \beta_l)_{\Omega_2}\quad  &
        K_{25} = (\Delta \psi_j, \mathbf{b} \cdot \Nabla \beta_l)_{\Omega_2}\\
        K_{26} = (\mathbf{b} \cdot \Nabla \psi_j, \mathbf{b} \cdot \Nabla \beta_l)_{\Omega_2}\quad  &
        \mathbf{F} = f_{i,j} \\
\end{array}
\end{equation}
\eirik{substitution}  into~\eqref{eqn:weak_form3c} yields the global system:
\one{
\begin{equation} \label{eqn:discrete_AD_1direction}
\begin{array}{c}
    (\kappa K_{11} \otimes K_{22} + \kappa K_{12} \otimes K_{21} + K_{12} \otimes K_{23} -\kappa \tau K_{13} \otimes K_{24} - \kappa \tau K_{12} \otimes K_{25} + \tau K_{12} \otimes K_{26})\mathbf{U} = \\ (K_{12} \otimes K_{22} + \tau K_{12} \otimes K_{24})\mathbf{F}.
\end{array}
\end{equation}
}
%Dirichlet boundary conditions are \eirik{enforced} by replacing the rows of the matrix corresponding to degrees of freedom on the boundary by the corresponding row of the identity matrix and the corresponding entry in the right hand side is replaced by the value of the Dirichlet condition at the boundary.

\subsection{Model Problem 4: Poisson Problem with Variable, Non-Separable Diffusivity}\label{sec:prob4}
\one{
To demonstrate that the methodology can be extended to more complicated settings where the problem is non-separable we will consider a problem similar to the one from Section~\ref{sec:model_prob} but now with a varying diffusion coefficient $\kappa$:
\begin{equation}
    \begin{split}
        -\nabla \cdot \kappa \nabla u = f, \quad \textrm{on} \quad \Omega \\
        u = u_D \quad \textrm{on} \partial \Omega,
    \end{split}
\end{equation}
\one{
where $\kappa$ is a bounded, continuous function of position $x$, but not of the solution variable $u$.  The derivation is essentially identical to Section~\ref{sec:model_prob} up to~\eqref{eqn:product_rule}. However, instead of~\eqref{eqn:product_rule} we have: }
\begin{equation}
           \int_{\Omega_1} \int_{\Omega_2}  \kappa \Nabla_1\phi_i \cdot \Nabla_1(\gamma_k)\,\psi_j\,\beta_l +  \phi_i\, \gamma_k\, \kappa \Nabla_2 \psi_j \cdot \Nabla_2\beta_l  \dy  \dx.
\end{equation}
The $\kappa$ term does not allow for full separability and the resulting Kronecker product structure as seen in previous cases can still construct a global system of equations via the following steps. In this case, the global stiffness matrix will be constructed corresponding to the following rearrangement of the above system:  
\begin{equation}\label{eq:problem4}
           \int_{\Omega_2} \psi_j\,\beta_l \int_{\Omega_1}   \kappa \Nabla_1\phi_i \cdot \Nabla_1\gamma_k\, \dx \dy +  \int_{\Omega_2}\Nabla_2 \psi_j \cdot \Nabla_2\beta_l \int_{\Omega_1} \kappa  \phi_i\, \gamma_k\,  \dx  \dy.
\end{equation}
Notice that the integrals over the domain $\Omega_1$ are only functions of $y \in {\Omega_2}$ since $\kappa$ is a function of both $x$ and $y$. For simplicity we can write:
\begin{equation}
           f_1(y) = \int_{\Omega_1}   \kappa(x,y) \Nabla_1\phi_i \cdot \Nabla_1\gamma_k\, \dx  \quad  f_2(y) = \int_{\Omega_1} \kappa (x,y) \phi_i\, \gamma_k\,  \dx.
\end{equation}
}

\one{To construct the global stiffness matrix, the only required task is the evaluation of  the integrands in~\eqref{eq:problem4}. This assembly procedure is not as efficient as the cases where the operators are completely separable between subdomains. First, let us denote the number of degrees of freedom in $\Omega_1 = N_1$ and the number of degrees of freedom in $\Omega_2 = N_2$. Furthermore, let us assume the quadrature rule being used only needs the function values at the degrees of freedom. Then, the algorithm can be summarized in the following steps: 
\begin{itemize}
    \item Compute $f_1$ and $f_2$ at all degrees pf freedom in the second subdomain $y_j \in \Omega_2$. This computation yields a set of $N_2$ sparse matrices of size $N_1 \times N_1$. Each sparse matrix represents the value of $f_1$,$f_2$ at a fixed point $y \in \Omega_2$.
    \item Evaluate each integral in~\eqref{eq:problem4} using the evaluations of $f_1$ and $f_2$ from the previous step. This  results in a global block structured matrix of dimension $N_1 N_2 \times N_1 N_2$ where each block will be a sparse $N_2 \times N_2$ matrix. 
    \item For efficiency, the $N_2 \times N_2$ blocks only need to be computed for the nonzero entries in each $N_1 \times N_1$ matrix.
    \item The right hand side is computed as in Section~\ref{sec:model_prob} with the Kronecker product.
    \item Modify global system to be consistent with boundary conditions where necessary.
\end{itemize}
}

\section{Numerical Verifications}\label{sec:numerical_results}
For each of the four problem \eirik{introduced in Section~\ref{sec:methods}, we consider and implement a specific test case} in FEniCS. \two{Since all of the above derivations only rely on integration of the subdomains, the implementation in FEniCS is possible without modification of the FEniCS codebase. Detailed tutorials for each of the following test cases are available on GitHub at~\url{https://github.com/Markloveland/FEniCS_Tensor_Product_Demos.git} in the form of Jupyter notebooks.} \eirik{To verify the developed framework, we investigate the $h-$convergence properties of the implemented methods by consideration of the  rate of convergence of the FE solutions. The test and trial spaces in all cases  consist of continuous Lagrange polynomials of degree 1.} \two{To do this for each test case, the $L^{\infty}$ and $L^2$ error norms are computed as the grids are uniformly refined. The $L^{\infty}$ is computed as:
\begin{equation}
    \|e\|_{L^{\infty}} = \max_{i,j\in N} |u_{exact}(x_i) - u_{i,j}|,
\end{equation}
 where $N$ is the set that contains the indices for all $i,j$ nodes, while $L^2$ error is computed as:
\begin{equation}
    \| e \|_{L^2} = \sqrt{\int_{\Omega} (u_{exact}-u)^2 dx}.
\end{equation}
The convergence rates between successive refinement steps are then computed as:
\begin{equation}
    \textrm{rate}_n  = \frac{ln(e_{n-1}/e_{n})}{ ln(d_{n-1}/d_{n})},
\end{equation}
where $n$ denotes the $n^{th}$ level of refinement and $d_n$ is the diameter of the element at the $n^{th}$ refinement level. For all of the finite element discretizations presented hereafter, we expect convergence rates to be close to 2 which is the optimal rate of convergence for Galerkin FE discretizations using linear polynomials, see, e.g.,Chapter 5 of the classical text by Carey and Oden~\cite{carey1983finite}.
}
%Each test case was run with several levels of grid refinement in order to investigate convergence rates. Errors were calculated with respect to the analytic solutions at each degree of freedom.

\subsection{Case 1: 4-D Poisson Equation }

\eirik{ We first consider the Poisson PDE in the high dimensional space of order four. We define the computational domain as } 
%Examining the problem (~\ref{eqn:model_prob1} ) and defining the domain as 
a tensor product between two unit squares, \eirik{i.e.,}
$\Omega = ((0,1)\times(0,1)) \times ((0,1) \times (0,1))$. %\eirik{In the 6-D case, we define the domain completely analogous as a product of two unit cubes.}
\eirik{We select} the forcing function defined as $f = \text{4}\pi^2 u_{exact}$  and the exact solution as $u_{exact} = \Pi^4_{i=1} sin(\pi x_i) $,
 where $x_i$ is a coordinate in the domain $\Omega$.  
 %\eirik{For simplicity in implementation in FEniCS, we consider homogeneous Dirichlet boundary conditions on} the  boundary $\partial \Omega$.  
 \eirik{In Table~\ref{tab:results_poisson}, the convergence data for the four  dimensional Poisson problem is listed. Note that the convergence of the FE solution is optimal, as  the rate of convergence of the root mean square error (RMSE) and $l_{\infty}$ norm approaches $\mathcal{O}(h^{p+1})$}.
%error was computed using $l_{\infty}$ and RMSE norms of the nodal values which are in table~\ref{tab:results_poisson}.
%For this case, the four dimensional Poisson problem, the RMSE converged at a rate on average of 1.56 with respect to h, while the $l_{\infty}$ converged at an average rate of 1.48...
%
%
\begin{table}[h!]
\centering
\caption{\label{tab:results_poisson}  Error estimation results for the 4D Poisson problem.}
\begin{tabular}{@{}llllll@{}}
\toprule
{dofs \hspace{6mm}} & { h \hspace{6mm}} & {$L_{\infty}$ \hspace{6mm}} & {$L_{\infty}$ rate} & {  $L^2$  \hspace{6mm}} & {$L^2$ rate}  \\
\midrule \midrule

256 &  0.333    & 1.40E-01  & -         & 5.62E-02     & - \\
1296 & 0.200    & 7.88E-02  & 1.13      & 2.32E-02     &  1.73\\
2401 & 0.167    & 6.55E-02  & 1.01      & 1.65E-02     & 1.88   \\
4096 & 0.143    & 4.50E-02  & 2.44      & 1.22-02      & 1.91 \\
6561 & 0.125    & 3.76E-02  & 1.34      & 9.47E-03     & 1.94. \\

%256 &  0.33333    & 0.140414487708  & 0.031603261 & - & - \\
%1296 & 0.20000    & 0.0788058171971 & 0.016107126 & 1.44 & 1.13 \\
%2401 & 0.16667 & 0.0655176797108 & 0.012103114 & 1.57 & 1.01  \\
%4096 & 0.14286 & 0.0449918036921 & 0.009390945 & 1.65 & 2.44 \\
%6561 & 0.12500 & 0.037599902 & 0.007483795 & 1.70 & 1.34\\
%1461 & 0.10000 & 0.02428415 & 0.005059920 & 1.75 & 1.96 \\
%20736 & 0.09091 & 0.019493556 & 0.004264705 & 1.79 & 2.31 \\
%28561 & 0.08333 & 0.016946887 & 0.003642222 & 1.81 & 1.61 \\

\bottomrule
\end{tabular}
\end{table}

\subsection{Case 2: 2D Space-Time Wave Equation}

\eirik{As a verification of the space-time wave model problem}~\eqref{eqn:model_prob2}, \eirik{we select} the space-time domain as the Cartesian product of a unit square \eirik{spatial} domain $\Omega$ and an interval time domain, i.e., $\Omega_T = ((0,1) \times (0,1)) \times (0,T)$.  The wave propagation speed is $c=1$ and \eirik{we consider a}  manufactured solution $u_{exact} = sin(x-ct) + sin(y-ct)$. \eirik{This solution is used to ascertain boundary and initial conditions needed to solve the resulting system of equations.}
%The boundary conditions were defined as $u(x,t) = u_{exact}$ on $\partial \Omega_1 \times [0,T]$, and the initial condition as $u(x,0) = sin(x) + sin(y)$. The boundary at $t=T$ was defined as 0 Neumann condition.
\eirik{In Table~\ref{tab:wave_spacetime}, the convergence results are presented along with the time interval element size, denoted by dt, and the space-time CFL number. The RMSE is observed to converge at the expected optimal rate, whereas the  }  $l_{\infty}$ error \eirik{exhibits a reduced rate for the finer meshes}. %The results are shown in Table~\ref{tab:wave_spacetime}.
\begin{table}[h!]
\centering
\caption{\label{tab:wave_spacetime}  Results for the Wave Space-Time problem.}
\begin{tabular}{@{}llllllll@{}}
\toprule
{dofs \hspace{6mm}} & { h \hspace{6mm}} & { dt \hspace{6mm}} & { CFL \hspace{6mm}} & {$L_{\infty}$ \hspace{6mm}} & {$L_{\infty}$ rate} & {$L^2$ error \hspace{15mm}} &{  $L^2$ rate }  \\
\midrule \midrule

200 &  0.250 & 0.14 & 0.57  & 2.25E-04  & -    & 7.80E-05 & -   \\
1215 & 0.125 & 0.07 & 0.57  & 5.06E-05  & 2.22 & 1.75E-05 & 2.15  \\
3718 & 0.083 & 0.05 & 0.57  & 2.53E-05  & 1.71 & 7.77E-06 & 2.01 \\
8381 & 0.063 & 0.04 & 0.57  & 1.31E-05  & 2.27 & 4.27E-06 &  2.08\\

%200 &  0.250 & 0.142857143 & 0.571428571  & 2.25E-04  & 5.83899E-05 & - & -   \\
%1215 & 0.125 & 0.071428571 & 0.571428571  & 5.06E-05  & 1.50656E-05 & 1.95 & 2.22  \\
%3718 & 0.083 & 0.047619048 & 0.571428571  & 2.53E-05  & 7.00547E-06 & 1.89 & 1.71 \\
%8381 & 0.063 & 0.035714286 & 0.571428571  & 1.31E-05  & 3.94971E-06 & 1.99 & 1.31 \\

\bottomrule
\end{tabular}
\end{table}

\subsection{Case 3: SUPG Stabilized Advection Dominated Advection Diffusion Equation}
\begin{figure}[h!]
    \centering
    \includegraphics[width=0.5\textwidth]{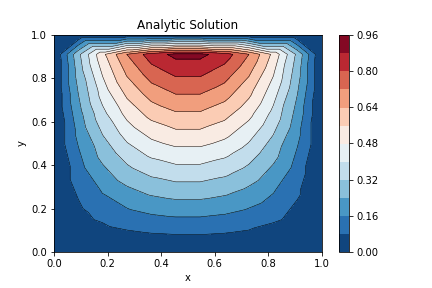}
    \caption{Analytic solution to Case 3 projected onto FE mesh with 289 dofs.}
    \label{fig:analytic_ADR}
\end{figure}

\begin{figure}[h!]
    \centering
    \includegraphics[width=0.5\textwidth]{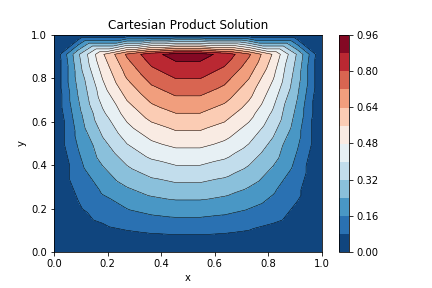}
    \caption{Cartesian product FE solution to Case 3 at 289 dofs.}
    \label{fig:cartesian_ADR}
\end{figure}

\eirik{As a final numerical verification, we consider a special case of advection dominated advection diffusion equation. In particular, we consider the case in which the advection 
acts in a single direction aligned with a coordinate axis, see~\eqref{eqn:discrete_AD_1direction}.    }
%Here we take the problem from~\eqref{eqn:model_problem3}   and use a test case as inspired by Egger and Sch{\"o}berl in~\cite{egger2010hybrid}. Hence, we use~\eqref{eq:fubinisplit} with a modified velocity term. 
 \eirik{We consider a case where the} domain $\Omega$ is a tensor product of 2 unit intervals: $\Omega = (0,1) \times (0,1)$, diffusivity constant $\kappa =\frac{1}{100}$, and \eirik{advection  vector is} $\mathbf{b} = (0,1)$. The analytic solution in this case \eirik{is inspired by the work of Egger and Sch{\"o}berl~\cite{egger2010hybrid}}:
\begin{equation} \label{eq:ADR_exact}
u_{exact} =  (-4(x-0.5)^2+1)\left[ y + \frac{e^{\frac{1}{\kappa} \cdot b_y \cdot y} - 1 }{1 - e^{\frac{1}{\kappa} \cdot b_y}}\right],
\end{equation}
which implies that the forcing term must be $f = -\kappa\,\Delta u_{exact} + \mathbf{b} \cdot \Nabla u_{exact}$ \eirik{and we enforce the corresponding} homogeneous Dirichlet on the boundary $\partial \Omega$.
\eirik{In Figures~\ref{fig:cartesian_ADR}, and~\ref{fig:analytic_ADR}, we show the approximate FE solution and the analytic solution, respectively. As expected, the SUPG stabilization results in a stable solution at this relatively coarse FE mesh.}
\eirik{In Table~\ref{tab:results_ADR}, the convergence results for this final case are presented. Both the RMSE and $l_{\infty}$ error converge at the expected optimal rates.} 
%The RMSE converged at an average rate of 1.72 while the $l_{\infty}$ converged at an average rate of 1.86. These rates are close to the expected optimal rate for finite element bases with polynomials of degree 1 which is of order 2. Also note, we started at higher refinement since convergence is not optimal unless boundary layer is resolved 
%
%
\begin{table}[h]
\centering
\caption{\label{tab:results_ADR}  Error estimation results for advection dominated advection diffusion problem.}
\begin{tabular}{@{}llllll@{}}
\toprule
{dofs \hspace{6mm}} & { h \hspace{6mm}} & {$L_{\infty}$ \hspace{6mm}} & {$L_{\infty}$ rate} & {$L^2$ \hspace{15mm}} & {  $L^2$ rate  } \\
\midrule \midrule

1089     & 0.03125    & 5.10E-04       & -    & 2.75E-04  & - \\
4225     & 0.015625   & 1.44E-04       & 1.83 & 7.75E-05  & 1.83 \\
4761     & 0.014703   & 1.28E-04       & 1.91 & 6.90E-05  & 1.92 \\
5329     & 0.013889   & 1.15E-04       & 1.91 & 6.18E-05  & 1.93  \\

%9    & 0.5     & 0.25            & 0.083333333 & -  & - \\
%16   & 0.333  & 0.118579208     & 0.046543091 & 1.44 & 1.84 \\
%25   & 0.25    & 0.073770219     & 0.03013949  & 1.51 & 1.65 \\
%36   & 0.2     & 0.047904246     & 0.02083341  & 1.65 & 1.93  \\
%49   & 0.167 & 0.033310117     & 0.01522303  & 1.73  & 1.99 \\
%256  & 0.067 & 0.0064096       & 0.00284602  & 1.83 & 1.80 \\
%676  & 0.04    & 0.002522659     & 0.001069127 & 1.92 & 1.83 \\
%2601 & 0.02   & 0.000649801     & 0.000274277 & 1.96 & 1.96 \\

\bottomrule
\end{tabular}
\end{table}
\subsection{Case 4: Poisson Equation with Variable Coefficient}
\one{We consider problem 4 from Section~\ref{sec:prob4} with the following set up: the domain is a Cartesian product of two unit intervals $\Omega = (0,1) \times (0,1)$. The coefficient $\kappa$ is a non-separable scalar function $\kappa = e^{\alpha xy}$ where $\alpha = 1$. The right hand side is set to $f=2\alpha^2(x^2+y^2)e^{2\alpha x y}$ and the analytic solution is $u = e^{\alpha(xy)}$. The boundary conditions are Dirichlet on the entire boundary and set to the exact solution. In Table~\ref{tab:results_poisson_vary} we present the corresponding convergence results and note that the convergence rates for both $L^{\infty}$ and $L^2$ errors are optimal at 2.}
\begin{table}[h!]
\centering
\caption{\label{tab:results_poisson_vary}  Error estimation results for the Poisson problem with a variable coefficient.}
\begin{tabular}{@{}llllll@{}}
\toprule
{dofs \hspace{6mm}} & { h \hspace{6mm}} & {$l_{\infty}$ \hspace{6mm}} & {$l_{\infty}$ rate} & {  $L^2$  \hspace{6mm}} & {$L^2$ rate}  \\
\midrule \midrule

25 &  0.25       & 7.34E-03  & -         & 3.96E-03     & - \\
81 & 0.125       & 1.72E-03  & 2.09      & 9.11E-04     &  2.12\\
289 & 0.0625     & 4.18E-04  & 2.04      & 2.23E-04     & 2.03   \\
1089 & 0.03125   & 1.04E-04  & 2.01      & 5.53-05      & 2.01 \\

\bottomrule
\end{tabular}
\end{table}

\section{Conclusions}\label{sec:Conclusions}

\eirik{In this paper, we have introduced and implemented tensor product FE routines  for high dimensional problems in FEniCS. This methodology allows us to extend the}
FEniCS library to domains with more than three dimensions so long as they are a Cartesian product of subdomains three or lower.
\eirik{To verify the developed methodology, we consider four test cases utilizing classical and stabilized FE methods. } For each test case, \eirik{we observe the expected } convergence to the analytic solutions with respect to grid refinement was demonstrated in both the $L^2$ and $L^{\infty}$ norms.

\eirik{We consider only linear PDEs here }  since this allowed for the explicit construction of a single linear system of algebraic equations. \eirik{Hence, future studies should investigate potential extensions to nonlinear PDEs.} Additionally, the global system of equations was solved naively by explicitly constructing the global stiffness matrix and inverting. \eirik{However,} the global matrix is sparse with highly structured blocks which should allow for faster solvers that would greatly reduce run time. \eirik{We refer to existing works~\cite{BIALECKI1993369,Gao_thesis,Bank1978}, where  related}  problems  were considered \eirik{and leave the consideration of such solvers for future studies.} 
\eirik{Further extensions to mixed FE methods, } such as those discussed in the book by Brezzi \emph{et al.}~\cite{brezzi2012mixed}. \eirik{Finally, } FE methods utilizing with discontinuous test/trial spaces \eirik{ could be considered due to their extensive use in engineering applications}. Furthermore, full integration of this method into the FEniCS API would be valuable for both simplicity of future implementations as well as for performance.

%There are several topics that could improve this research in the future. An interesting topic that has garnered some attention is creating a fast solver that can utilize the highly structured nature of the global linear systems.
% Another improvement to this research would be to extend this tensor product method to work for nonlinear PDEs and to extend the work to different mixed finite element methods

\section{Acknowledgements}

Author Loveland has been supported by the CSEM Fellowship from the Oden Institute at the University of Texas at Austin. Authors Loveland, Valseth, and Dawson have been supported by the United States National Science Foundation - NSF PREEVENTS Track
2 Program, under NSF Grant Number 1855047 and the Department of Homeland Security Coastal Resilience Center research project "Accurate and Fast Wave Modeling and Coupling with ADCIRC".
Author Lukac has been supported by the University of Oregon. The authors would also like to thank the reviewers of this manuscript for their time and thoughtful suggestions.

\bibliographystyle{elsarticle-num}
 \bibliography{references}

\begin{thebibliography}{10}
\expandafter\ifx\csname url\endcsname\relax
  \def\url#1{\texttt{#1}}\fi
\expandafter\ifx\csname urlprefix\endcsname\relax\def\urlprefix{URL }\fi
\expandafter\ifx\csname href\endcsname\relax
  \def\href#1#2{#2} \def\path#1{#1}\fi

\bibitem{rathgeber2016firedrake}
F.~Rathgeber, D.~A. Ham, L.~Mitchell, M.~Lange, F.~Luporini, A.~T. McRae, G.-T.
  Bercea, G.~R. Markall, P.~H. Kelly, Firedrake: automating the finite element
  method by composing abstractions, ACM Transactions on Mathematical Software
  (TOMS) 43~(3) (2016) 1--27.

\bibitem{bangerth2007deal}
W.~Bangerth, R.~Hartmann, G.~Kanschat, deal. ii—a general-purpose
  object-oriented finite element library, ACM Transactions on Mathematical
  Software (TOMS) 33~(4) (2007) 24--es.

\bibitem{anderson2021mfem}
R.~Anderson, J.~Andrej, A.~Barker, J.~Bramwell, J.-S. Camier, J.~Cerveny,
  V.~Dobrev, Y.~Dudouit, A.~Fisher, T.~Kolev, et~al., Mfem: A modular finite
  element methods library, Computers \& Mathematics with Applications 81 (2021)
  42--74.

\bibitem{alnaes2015fenics}
M.~Aln{\ae}s, J.~Blechta, J.~Hake, A.~Johansson, B.~Kehlet, A.~Logg,
  C.~Richardson, J.~Ring, M.~E. Rognes, G.~N. Wells, The fenics project version
  1.5, Archive of Numerical Software 3~(100) (2015).

\bibitem{renard2020getfem}
Y.~Renard, K.~Poulios, Getfem: Automated fe modeling of multiphysics problems
  based on a generic weak form language, ACM Transactions on Mathematical
  Software (TOMS) 47~(1) (2020) 1--31.

\bibitem{brenner2008mathematical}
S.~C. Brenner, L.~R. Scott, L.~R. Scott, The mathematical theory of finite
  element methods, Vol.~3, Springer, 2008.

\bibitem{quarteroni2010numerical}
A.~Quarteroni, R.~Sacco, F.~Saleri, Numerical mathematics, Vol.~37, Springer
  Science \& Business Media, 2010.

\bibitem{ern2013theory}
A.~Ern, J.-L. Guermond, Theory and practice of finite elements, Vol. 159,
  Springer Science \& Business Media, 2013.

\bibitem{Bank1978}
R.~E. Bank, Efficient algorithms for solving tensor product finite element
  equations, Numerische Mathematik 31~(1) (1978) 49--61.
\newblock \href {https://doi.org/10.1007/BF01396013}
  {\path{doi:10.1007/BF01396013}}.

\bibitem{BAKER1981215}
A.~Baker, M.~Soliman, On the accuracy and efficiency of a finite element tensor
  product algorithm for fluid dynamics applications, Computer Methods in
  Applied Mechanics and Engineering 27~(2) (1981) 215--237.
\newblock \href {https://doi.org/https://doi.org/10.1016/0045-7825(81)90150-X}
  {\path{doi:https://doi.org/10.1016/0045-7825(81)90150-X}}.

\bibitem{DU2013181}
K.~Du, W.~Sun, X.~Zhang, Arbitrary high-order c0 tensor product {G}alerkin
  finite element methods for the electromagnetic scattering from a large
  cavity, Journal of Computational Physics 242 (2013) 181--195.
\newblock \href {https://doi.org/https://doi.org/10.1016/j.jcp.2013.02.015}
  {\path{doi:https://doi.org/10.1016/j.jcp.2013.02.015}}.

\bibitem{mcrae2016automated}
A.~T. McRae, G.-T. Bercea, L.~Mitchell, D.~A. Ham, C.~J. Cotter, Automated
  generation and symbolic manipulation of tensor product finite elements, SIAM
  Journal on Scientific Computing 38~(5) (2016) S25--S47.

\bibitem{munch2020hyperdeal}
P.~Munch, K.~Kormann, M.~Kronbichler, hyper.deal: An efficient, matrix-free
  finite-element library for high-dimensional partial differential equations
  (2021).

\bibitem{becker1981finite}
E.~B. Becker, G.~F. Carey, J.~T. Oden, Finite elements: an introduction,
  Vol.~1, Prentice Hall, 1981.

\bibitem{Loscher2021}
R.~Löscher, O.~Steinbach, M.~Zank, Numerical results for an unconditionally
  stable space-time finite element method for the wave equation (2021).
\newblock \href {https://doi.org/10.48550/ARXIV.2103.04324}
  {\path{doi:10.48550/ARXIV.2103.04324}}.

\bibitem{brooks1982streamline}
A.~N. Brooks, T.~J. Hughes, Streamline upwind/{P}etrov-{G}alerkin formulations
  for convection dominated flows with particular emphasis on the incompressible
  navier-stokes equations, Computer methods in applied mechanics and
  engineering 32~(1-3) (1982) 199--259.

\bibitem{carey1983finite}
G.~F. Carey, J.~T. Oden, Finite Elements: A Second Course; Graham F. Carey and
  J. Tinsley Oden, Prentice-hall, 1983.

\bibitem{egger2010hybrid}
H.~Egger, J.~Sch{\"o}berl, A hybrid mixed discontinuous {G}alerkin
  finite-element method for convection--diffusion problems, IMA Journal of
  Numerical Analysis 30~(4) (2010) 1206--1234.

\bibitem{BIALECKI1993369}
B.~Bialecki, G.~Fairweather, Matrix decomposition algorithms for separable
  elliptic boundary value problems in two space dimensions, Journal of
  Computational and Applied Mathematics 46~(3) (1993) 369--386.
\newblock \href {https://doi.org/https://doi.org/10.1016/0377-0427(93)90033-8}
  {\path{doi:https://doi.org/10.1016/0377-0427(93)90033-8}}.

\bibitem{Gao_thesis}
L.~Gao, \href{http://hdl.handle.net/10754/303766}{Kronecker products on
  preconditioning} (2013).
\newblock \href {https://doi.org/10.25781/KAUST-8S7R9}
  {\path{doi:10.25781/KAUST-8S7R9}}.
\newline\urlprefix\url{http://hdl.handle.net/10754/303766}

\bibitem{brezzi2012mixed}
F.~Brezzi, M.~Fortin, Mixed and hybrid finite element methods, Vol.~15,
  Springer Science \& Business Media, 2012.

\end{thebibliography}
\end{document}